\documentclass{article}
\usepackage{amssymb}

\usepackage{graphicx}
\usepackage{amsmath}

\newtheorem{theorem}{Theorem}

\newtheorem{definition}[theorem]{Definition}

\newtheorem{lemma}[theorem]{Lemma}

\input{tcilatex}

\begin{document}

\begin{center}
\textbf{NOMBRE\ DE\ COMPOSANTES\ CONNEXES\ }

\textbf{D'UNE\ VARIETE\ REELLE}$\;$\textbf{ET \ }$\mathbb{R}-\mathbf{PLACES}$%
\textbf{\ }

\bigskip

Danielle GONDARD-COZETTE

danielle.gondard@imj-prg.fr

\bigskip 
\end{center}

\textbf{Abstract. }

\textit{The purpose of this paper is to present results and open problems
related to }$R-$\textit{places. The first section recalls basic
facts, the second introduces }$R-$\textit{places\ and their
relationship with orderings and valuations. }

\textit{The third part involves Real Algebraic Geometry and gives results
proved using the space of \ }$R-$\textit{places. Theorem 14 gives
explicitly, in terms of the function field of the variety, the number of
connected components of a non-empty smooth projective real variety.}

\textit{The fourth and fifth parts are devoted to the links with the real
holomorphy rings and the valuation fans. Then we present an approach to
abstract real places and conclude with some open questions.}

\bigskip

-----

\bigskip

\textbf{R\'{e}sum\'{e}.}

Le but de cet article est de rassembler des r\'{e}sultats et des questions
concernant un objet mal connu mais souvent utile, les $R-$places$,$
et de fournir les r\'{e}f\'{e}rences correspondantes$.$

Nous commen\c{c}ons par rappeler les notions de base, puis pr\'{e}sentons
quelques r\'{e}sultats de g\'{e}om\'{e}trie r\'{e}elle obtenus gr\^{a}ce aux 
$R-$places. En particulier le th\'{e}or\`{e}me 14 donne
explicitement, en termes de corps des fonctions de la vari\'{e}t\'{e}, le
nombre de composantes connexes d'une vari\'{e}t\'{e} r\'{e}elle projective
lisse et non vide.

Nous montrons ensuite quels sont les liens entre les $R-$places et
d'autres objets comme les ordres de niveau sup\'{e}rieur, les \'{e}ventails
valu\'{e}s et l'anneau d'holomorphie r\'{e}el. Enfin nous passons au cas des 
$R-$places dans des espaces d'ordres abstraits et concluons par
quelques questions ouvertes.

\section{Rappels}

\subsection{\protect\bigskip Espaces des ordres d'un corps}

Rappelons le r\'{e}sultat bien connu obtenu par Artin-Schreier en 1927%
{\large \ : }

\begin{theorem}
Un corps $K$ commutatif ordonnable est caracteris\'{e} par $-1\notin \sum
K^{2},$ o\`{u} $\sum K^{2}$ est l'ensemble des sommes finies de carr\'{e}s
d'\'{e}l\'{e}ments de $K$.
\end{theorem}

Cet ensemble $\sum K^{2}$ est aussi \'{e}gal \`{a} l'ensemble form\'{e} des
\'{e}l\'{e}ments positifs dans tous les ordres du corps.

On sait qu'il existe alors $P\subset K$ tel que $P+P\subseteq P$, $%
P.P\subset P$, $-1\notin P$, $P\cup -P=K$ ; $P$ est un ordre total
compatible avec la structure de corps de $K.$

\begin{definition}
On d\'{e}signe par $\chi (K)$ l'espace des ordres $P$ de $K$. La topologie
usuelle sur $\chi (K)$ est celle de Harrison engendr\'{e}e par les
ouverts-ferm\'{e}s 
\begin{equation*}
H(a)=\{P\in \chi (K)\mid a\in P\}.
\end{equation*}
\end{definition}

\QTP{Body Math}
Muni de cette topologie $\chi $ $(K)$ est un espace compact, totalement
disconnexe.

\QTP{Body Math}
Il a \'{e}t\'{e} montr\'{e} par Craven [8] que tout espace compact
totalement disconnexe \'{e}tait hom\'{e}omorphe \`{a} l'espace $\chi (K)$
des ordres d'un corps $K.$

\subsection{Les places r\'{e}elles}

\begin{definition}
Une place sur un corps $K$ est une application 
\begin{equation*}
\varphi :K\rightarrow F\cup \left\{ \infty \right\}
\end{equation*}
o\`{u} $F$ est un corps, telle que $\varphi (1)=1$ et qui satisfait les
r\`{e}gles d'homorphisme usuelles pour la somme et pour le produit.
\end{definition}

\begin{definition}
On appelle place r\'{e}elle, une place telle que le corps $F$ est ordonnable
(ou r\'{e}el-clos). Si $F=$ $\mathbb{R}$, la place r\'{e}elle est
appel\'{e}e une $\mathbb{R-}$place.
\end{definition}

Toutes les $\mathbb{R-}$places d'un corps peuvent \^{e}tre obtenues \`{a}
partir de l'espace des ordres du corps $\chi (K)$ en utilisant certaines
valuations r\'{e}elles.

\subsection{Les valuations r\'{e}elles}

\begin{definition}
Une valuation, au sens de Krull, est une application surjective $v:K^{\ast
}\twoheadrightarrow \Gamma ,$ $o\grave{u}$ $\Gamma =(\Gamma ,+,\leq )$ $est$ 
$un$ $groupe$ $abelien$ $totalement$ $ordonn\acute{e},$telle que

$\forall x,y\in K^{\ast },$ $v(xy)=v(x)+v(y)$.

$\forall x,y\in K^{\ast }$, $x+y\neq 0,$ $v(x+y)\geq \min (v(x),v(y)).$
\end{definition}

L'anneau de la valuation $v$ est d\'{e}fini par : 
\begin{equation*}
A=\{a\in K\mid a=0\text{ }ou\ v(a)\geqslant 0\},
\end{equation*}
et l'id\'{e}al maximal de $A$ est donn\'{e} par : 
\begin{equation*}
\QTR{sl}{m}=\{a\in K\mid a=0\text{ }ou\ v(a)>0\};
\end{equation*}
le groupe des unit\'{e}s est 
\begin{equation*}
A^{\mathrm{\ast }}=A\backslash \QTR{sl}{m,}
\end{equation*}
et le corps r\'{e}siduel 
\begin{equation*}
k_{v}=A/\QTR{sl}{m.}
\end{equation*}

\begin{definition}
La valuation $v$ est \textit{r\'{e}elle }si et seulement si le corps
r\'{e}siduel $k$ est ordonnable.
\end{definition}

Un corps commutatif $K$ admet des valuations r\'{e}elles si et seulement si
il est ordonnable (on dit aussi formellement r\'{e}el).

Plus pr\'{e}cis\'{e}ment, si $P$ est un ordre donn\'{e} de $K$, alors
l'enveloppe convexe de $\mathbb{Q}$ dans $K$ 
\begin{equation*}
A(P)=\{a\in K\mid \exists r\in \mathbb{Q},-r\leq _{P}a\leq _{P}r\}
\end{equation*}
est un anneau de valuation et 
\begin{equation*}
I(P)=\{a\in K\mid \forall r\in \mathbb{Q}^{\ast },-r\leq _{P}a\leq _{P}r\}
\end{equation*}
est son id\'{e}al maximal ; $P$ induit sur le corps r\'{e}siduel $\
k=A(P)/I(P)$ un ordre archim\'{e}dien $\overline{P}.$

\bigskip

\section{Les $\mathbb{R-}$places}

\subsection{$\mathbb{R}-$place associ\'{e}e \`{a} un ordre}

Pour une pr\'{e}sentation tr\`{e}s d\'{e}taill\'{e}e de ces notions on
pourra consulter [14], et pour d'autres r\'{e}sultats [17].

Soit $K$ un corps ordonnable, $P$ un ordre de $K$ $;$ on sait d'apr\`{e}s ce
qui pr\'{e}c\`{e}de que $(k_{v,}\overline{P})$ se plonge avec unicit\'{e}
dans $(\mathbb{R},\mathbb{R}^{2})$ ; on note $i$ ce plongement et $\pi $
l'application canonique de $K$ dans $k_{v}\cup \{\infty \}$ (o\`{u} si $%
a\notin A(P),$ alors $\pi (a)=\infty ).$

\begin{definition}
La $\mathbb{R}-place$ associ\'{e}e \`{a} $P$ est $\lambda _{P}:K\rightarrow 
\mathbb{R}\cup \{\infty \}$ d\'{e}finie par le diagramme commutatif :
\end{definition}

$\qquad \qquad \qquad \qquad \qquad \qquad K\qquad \overset{\lambda _{P}}{%
\longrightarrow }\qquad \mathbb{R}\cup \{\infty \}$ \ \ \ 

\bigskip

$\qquad \qquad \qquad \qquad \qquad \qquad \pi \searrow \qquad \qquad
\nearrow i$

\bigskip \qquad \qquad \qquad \qquad \qquad \qquad \qquad $k_{v}\cup
\{\infty \}$

\bigskip

\QTP{Body Math}
Explicitement : $\lambda _{P}(a)=\infty $ si $a\notin A(P),$ et si $a\in
A(P) $, alors 
\begin{equation*}
\lambda _{P}(a)=\inf \{r\in \mathbb{Q}\mid a\leq _{P}r\}=\sup \{r\in \mathbb{%
Q}\mid r\leq _{P}a\}.
\end{equation*}

\subsection{L'espace des $\mathbb{R}-$places}

\begin{definition}
L'espace des $\mathbb{R}$-places d'un corps est $M(K)=\{\lambda _{P}\mid
P\in \chi (K)\},$ o\`{u} $\chi (K)$ d\'{e}signe l'espace des ordres du corps
K.
\end{definition}

\QTP{Body Math}
L'espace $M(K)$ est muni de la topologie la plus grossi\`{e}re rendant
continues les applications \'{e}valuations d\'{e}finies \ pour tout $a\in K$
par 
\begin{equation*}
e_{a}\text{{\large \ }}{\large :}\text{ }M(K)\longrightarrow \mathbb{R}\cup
\{\infty \}
\end{equation*}
\begin{equation*}
\lambda \mapsto \lambda (a)
\end{equation*}

\QTP{Body Math}
$M(K)$ muni de cette topologie est un espace compact s\'{e}par\'{e}, et
l'application 
\begin{equation*}
\Lambda :\chi (K){\large \ }\longrightarrow M(K)
\end{equation*}
\begin{equation*}
P\mapsto \lambda _{P}
\end{equation*}
est continue, ferm\'{e}e et surjective.

\QTP{Body Math}
La topologie de $M(K)$ est aussi la topo$\log $ie quotient h\'{e}rit\'{e}e
de $\chi (K)${\large .}

\subsection{Lien avec les ordres de niveau sup\'{e}rieur (Becker [3])}

\begin{definition}
\textit{\ Soit }$K$ un corps ordonnable commutatif. $P\subset K$ est un\emph{%
\ }ordre de niveau exact $n$ si $0,1\in P,$ $-1\notin P,$ $P+P\subset P,$ $%
P.P\subset P,$ $P^{\ast }$ est un sous-groupe de $K^{\ast },$ et $K^{\ast
}/P^{\ast }\simeq \mathbb{Z}/2n\mathbb{Z}$.
\end{definition}

Les ordres de niveau 1 sont les ordres totaux usuels.

Les ordres de niveau sup\'{e}rieur peuvent aussi \^{e}tre compris en
utilisant des signatures : $P=\ker \sigma \cup \{0\},$ o\`{u} $\sigma
:K^{\ast }\rightarrow \mu _{2n}$ est un morphisme de groupes ab\'{e}liens de
noyau additivement ferm\'{e}, et o\`{u} $\mu _{2n}$ d\'{e}signe l'ensemble
des racines $2n$-i\`{e}mes de l'unit\'{e}.

\bigskip

A tout ordre de niveau sup\'{e}rieur on associe, comme pour un ordre usuel,
une unique $\mathbb{R-}$place ; un tel ordre \'{e}tant un cas particulier
d'\'{e}ventail valu\'{e} nous renvoyons le lecteur au \S 5\ plus loin.

\bigskip

Ces ordres de niveau sup\'{e}rieur pr\'{e}sentent des liens importants avec
les sommes de puissances.

Dans toute la suite $\sum K^{2p}$ d\'{e}signe l'ensemble des sommes finies
de puissances $2p$-i\`{e}mes d'\'{e}l\'{e}ments de $K.$

\begin{theorem}
Sont \'{e}quivalents ($p$\ premier) :

(1) $\sum K^{2}\neq \sum K^{2p}$

(2) $K$ admet un ordre de niveau exact $p$.
\end{theorem}

\begin{theorem}
Pour tout entier $n\geq 1$ $\ $on a 
\begin{equation*}
\sum K^{2n}=\underset{\text{{\small \ }}P\text{{\small \ ordre dont le
niveau divise }}n}{\cap }P
\end{equation*}
\end{theorem}

\textbf{Exemple : }Si $K=\mathbb{R((X))},$ les deux ordres usuels sont 
\begin{equation*}
P_{+}=K^{2}\cup XK^{2}\ {\small et\ }P_{-}=K^{2}\cup -XK^{2}
\end{equation*}
et pour tout\textit{\ }$p$\textit{\ }premier il y a deux ordres de niveau%
\textit{\ }$p\ $\textit{:} 
\begin{equation*}
P_{p,+}=K^{2p}\cup X^{p}K^{2p}{\small \ et\ }P_{p,-}=K^{2p}\cup -X^{p}K^{2p}
\end{equation*}
Tous ces ordres sont associ\'{e}s \`{a} l'unique $\mathbb{R-}$place de $%
\mathbb{R}((X))${\small , }et pour la valuation associ\'{e}e, ils induisent
tous sur le corps r\'{e}siduel le m\^{e}me ordre archim\'{e}dien.

\bigskip

Il r\'{e}sulte des travaux de Becker le th\'{e}or\`{e}me suivant :

\begin{theorem}
Les propositions suivantes sont \'{e}quivalentes \textit{:}

(1) $\forall a\in K$\ \ $a^{2}\in \sum K^{4}$\ ;

(2) toute valuation r\'{e}elle de $K$\ a un groupe des valeurs 2-divisible ;

(3) $K\ $n'admet pas d'ordre de niveau exact $2$.
\end{theorem}

En corollaire nous obtenons que $\lambda _{P}=\lambda _{Q}$ si et seulement
si $P$\ et $Q$\ sont le d\'{e}but d'une cha\^{i}ne 2-primaire d'ordres de
niveau sup\'{e}rieur (une telle cha\^{i}ne a \'{e}t\'{e} d\'{e}finie par
Harman [12] comme ($P_{n})=(P_{0,}P_{1},...,P_{n},...)$, $P_{0\text{ }}$%
ordre usuel, $P_{n}$\ ordre de niveau exact $2^{n-1}$\ tel que $P_{n}\cup
-P_{n}=(P_{0}\cap P_{n-1})\cup -(P_{0}\cap P_{n-1}))$

\QTP{Body Math}
$\bigskip $

\QTP{Body Math}
L'application $\Lambda $\ $:\chi (K)\ \longrightarrow M(K)$ est donc une
bijection si et seulement si le corps $K$ n'admet aucun ordre de niveau
exact 2.

\section{Une utilisation des $\mathbb{R}-$places\ en g\'{e}om\'{e}trie
r\'{e}elle}

\subsection{{\protect\large Un crit\`{e}re de s\'{e}paration des composantes
connexes dans }$M(K)$}

\begin{theorem}
(Becker-Gondard [6]) : Les $\mathbb{R-}places$ $\lambda _{P}$ et $\lambda
_{Q}$ sont dans deux composantes connexes distinctes de $M(K)$ si et
seulement si : 
\begin{equation*}
\exists b\in K^{\ast }\ (b\in P\cap -Q\ \text{et}\ b^{2}\in \sum K^{4}).
\end{equation*}
\end{theorem}

\textbf{Preuve} :

Ce crit\`{e}re est obtenu par la th\'{e}orie des ordres de niveau
sup\'{e}rieur, plus pr\'{e}cis\'{e}ment des ordres de niveau $2$.

On a $\chi (K)=H(a)\cup H(-a)$ et $H(a)\cap H(-a)=\varnothing ,$ mais $%
\Lambda (H(a))\cap \Lambda (H(-a))$ peut ne pas \^{e}tre vide.

Cependant, s'il existe $b\notin \sum K^{2}$ et tel que $b^{2}\in \sum K^{4},$
alors il n'existe pas $P\in H(b)$ et $Q\in H(-b)$ tels que $\lambda
_{P}=\lambda _{Q}$ .

Sinon $b\notin (P\cap Q)\cup -(P\cap Q)$ et $\lambda _{P}=\lambda _{Q}$
impliquent, comme il a \'{e}t\'{e} dit en fin de \S\ 2, qu'il existe un
ordre de niveau $2,$ $P_{2},$ tel que 
\begin{equation*}
P_{2}\cup -P_{2}=(P\cap Q)\cup -(P\cap Q)
\end{equation*}
avec $b\notin P_{2}\cup -P_{2}$, d'o\`{u} $b^{2}\notin P_{2}$ $,$ et donc $%
b^{2}\notin \sum K^{4}=\cap P_{2,i},$ o\`{u} $P_{2,i}$ parcourt l'ensemble
des ordres dont le niveau divise $2.$

\bigskip

Supposons alors que $\lambda _{P}$ et $\lambda _{Q}$ sont dans la m\^{e}me
composante connexe $C$ de $M(K)$ (avec $P\neq Q),$ et qu'il existe $b\in
P\cap -Q$ avec $b^{2}\in \sum K^{4}$ $;$ $\Lambda $ \'{e}tant ferm\'{e}e $%
C\cap \Lambda (H(b)),$ et $C\cap \Lambda (H(-b))$ forment une partition de $%
C $ en deux ferm\'{e}s non vides, impossible.

\bigskip

R\'{e}ciproquement :

Si $\lambda _{P}$ et $\lambda _{Q}$ sont dans $C$ et $C^{\prime },$ deux
composantes connexes distinctes de $M(K)$, $M(K)$ \'{e}$\tan $t compact
s\'{e}par\'{e} il existe un ouvert-ferm\'{e} $U$ $\supset C$ et $%
U^{c}=M(K)\backslash U\supset C^{\prime }$.

Soient $X=\Lambda ^{-1}(U)$ et $Y=\Lambda ^{-1}(U^{c})$ $;$ $X$ et $Y$
forment une partition de $\chi (K)$ $;$ $\Lambda $ \'{e}tant surjective on a
: 
\begin{equation*}
\Lambda ^{-1}(\Lambda (\Lambda ^{-1}(U)))=\Lambda ^{-1}(U)
\end{equation*}
donc $\Lambda ^{-1}(\Lambda (X))=X,$ et de m\^{e}me $\Lambda ^{-1}(\Lambda
(Y))=Y.$

Le lemme ci-apr\`{e}s de Harman donne alors l'existence de $b$ tel que $%
X=H(b)$ et $Y=H(-b)$ avec\ $b^{2}\in \sum K^{4}$, donc on a $b\in P\cap -Q$
avec $b^{2}\in \sum K^{4}.$

\bigskip

\QTP{Body Math}
\textbf{Lemme de Harman} ([12]) : Si $\chi (K)=\chi _{1}\cup \chi _{2}$,
o\`{u} $\chi _{1}$ et $\chi _{2}$ sont des ouverts-ferm\'{e}s disjoints tels
que $\Lambda ^{-1}(\Lambda (\chi _{1}))=\chi _{1}$ $et$ $\Lambda
^{-1}(\Lambda (\chi _{2}))=\chi _{2},$ alors il existe $a$ tel que $\chi
_{1}=H(a)$ et $\chi _{2}=H(-a).$

\subsection{Nombre de composantes connexes d'une vari\'{e}t\'{e} r\'{e}elle}

\begin{theorem}
(Becker-Gondard[6]) : Soit $Y$ une vari\'{e}t\'{e} projective lisse non vide
sur $\mathbb{R}$ , de corps des fonctions $K=$ $\mathbb{R(}Y\mathbb{)}.$
Alors $\mid \pi _{0}(Y(\mathbb{R}))\mid ,$ le nombre de composantes connexes
de $Y(\mathbb{R)},$ est donn\'{e} par : 
\begin{equation*}
\mid \pi _{0}(Y(\mathbb{R}))\mid =1+\log _{2}[(K^{\ast 2}\cap \sum
K^{4}):(\sum K^{\ast 2})^{2}]
\end{equation*}
\end{theorem}

Remarques :

Ce r\'{e}sultat est dans l'esprit de celui de Harnack qui majore le nombre
de composates connexes d'une courbe projective lisse $V(\mathbb{R)}$ par $%
g+1 $, o\`{u} $g$ est le genre de V ; mais ici nous avons une formule avec
\'{e}galit\'{e}. Le th\'{e}or\`{e}me met bien en \'{e}vidence le fait connu
que le nombre de composantes connexes est un invariant birationnel parmi les
vari\'{e}t\'{e}s lisses.

\bigskip

La premi\`{e}re preuve de ce r\'{e}sultat peut \^{e}tre trouv\'{e}e dans [6].

\bigskip

Deux nouvelles preuves de ce th\'{e}or\`{e}me ont \'{e}t\'{e} trouv\'{e}es
en 2003-2004 par Jean-Louis Colliot-Th\'{e}l\`{e}ne [7] et par Claus
Scheiderer [16].

\bigskip

Dans la preuve originelle, le th\'{e}or\`{e}me r\'{e}sulte des deux lemmes
ci-dessous qui utilisent les composantes connexes de l'espace des $\mathbb{R}%
-$places\ $M(K).$

\begin{lemma}
Soit $Y$ une vari\'{e}t\'{e} projective lisse sur $\mathbb{R}$, de corps des
fonctions $K=$ $\mathbb{R}(Y).$ Alors $\mid \pi _{0}(Y(\mathbb{R}))\mid $ le
nombre de composantes connexes de $Y$ est \'{e}gal \`{a} : 
\begin{equation*}
\mid \pi _{0}(Y(\mathbb{R}))\mid =\mid \pi _{0}(M(\mathbb{R(}Y)))\mid .
\end{equation*}
\end{lemma}

\begin{lemma}
Pour tout corps ordonnable $K$ : 
\begin{equation*}
\mid \pi _{0}(M(K))\mid =1+\log _{2}[(K^{\ast 2}\cap \sum K^{4}):(\sum
K^{\ast 2})^{2}].
\end{equation*}
\end{lemma}

\QTP{Body Math}
\textbf{Esquisse de preuve du lemme 15.}

On utilise l'application centre $c:$ $M(K)\rightarrow Y($ $\mathbb{R)},$
d\'{e}finie par $x=c(\lambda )=c(V_{\lambda })$ l'unique point ($Y$
projective) dont l'anneau local $\frak{o}_{x}$ est domin\'{e} par $%
V_{\lambda },$ l'anneau de valuation associ\'{e}e \`{a} la $\mathbb{R}-$%
place $\lambda $.

- Dans ce cas il est connu [par ex. Bochnak-Coste-Roy, G\'{e}om\'{e}trie Alg%
\'{e}brique R\'{e}elle, Prop. 7.6.2 (ii), p. 133] que $c$ est surjective,
les points centraux \'{e}tant l'adh\'{e}rence des points r\'{e}guliers. Et
on peut montrer que $c$ est continue.

- Br\"{o}cker a montr\'{e} (non publi\'{e}) que la fibre d'un point central
a un nombre fini de composantes connexes, et que si $x$ est r\'{e}gulier
alors elle est connexe.

Enfin on utilise le lemme suivant : si une application entre deux espaces
compacts $X$ et $Y$ est continue surjective et que chaque fibre est connexe
alors elle induit une bijection entre $\pi _{0}(X)$ et $\pi _{0}(Y)$.

\bigskip

\QTP{Body Math}
\textbf{Esquisse de preuve du lemme 16.}

On montre que $\mid \pi _{0}(M(K))\mid =\log _{2}[E:E^{+}]$ o\`{u} $E$ est
le groupe des unit\'{e}s de l'anneau d'holomorphie r\'{e}el $H(K)$ et $%
E^{+}=E\cap \sum K^{2}.$

Ensuite, on peut prouver que le groupe quotient ($K^{\ast 2}\cap \sum
K^{4})/(\sum K^{\ast 2})^{2}$ est isomorphe \`{a} $E/(E^{+}\cup -E^{+}).$

\section{Les $\mathbb{R}-$places et l'anneau d'holomorphie r\'{e}el}

\begin{definition}
On appelle anneau d'holomorphie r\'{e}el,\emph{\ }et\emph{\textit{\ }}on note%
\emph{\ }$H(K),$ l'anneau intersection de tous les anneaux de valuation
r\'{e}elle sur $K.$

On a aussi $H(K)=\underset{P\in \chi (K)}{\cap }A(P),$ et 
\begin{equation*}
H(K)=A(\sum K^{2})=\{a\in K\mid \exists n\in \mathbb{N},n\geq 1\text{ }tel%
\text{ }que\text{ }n\pm a\in \sum K^{2}\}.
\end{equation*}
\end{definition}

$H(K)$ est un anneau de Pr\"{u}fer (anneau $R\subset K$ tel que pour tout
id\'{e}al premier $p$ le localis\'{e} $R_{p}$ est un anneau de valuation de $%
K$), de corps des quotients $K.$

On notera dans la suite

\begin{center}
$Sper(H(K))=\{\alpha =(\QTR{sl}{p,}\overline{\QTR{sl}{\alpha }}),$ $\QTR{sl}{%
p\in specH(K),}$ $\overline{\QTR{sl}{\alpha }}$ ordre de $quot(H(K)/\QTR{sl}{%
p)\}}$
\end{center}

\QTP{Body Math}
le spectre r\'{e}el de l'anneau d'holomorphie r\'{e}el de $K$.

\begin{theorem}
(Becker-Gondard [6]) : On a le diagramme commutatif suivant :
\end{theorem}

\begin{center}
$\chi (K)$ $\qquad \ \overset{speri}{\longrightarrow }$ $\qquad \
MinSperH(K) $

$\downarrow \Lambda \qquad $ $\qquad \qquad \qquad \qquad \qquad \qquad
\downarrow sp$

$M(K)\overset{res}{\rightarrow }Hom(H(K),\mathbb{R})\overset{j}{\rightarrow }%
MaxSperH(K)$
\end{center}

\QTP{Body Math}
\textit{o\`{u} les applications horizontales sont des hom\'{e}omorphismes,
et les verticales des surjections continues.}

\bigskip

Les applications du diagramme ci-dessus sont d\'{e}finies comme suit :

\bigskip

$\Lambda :\chi (K)\longrightarrow M(K)$ est donn\'{e}e par $P\mapsto \lambda
_{P}$ ;

\bigskip

$speri:\chi (K)\longrightarrow MinSperH(K)$ est donn\'{e}e par $P\mapsto
P\cap H(K)$ ;

\bigskip

$sp:MinSperH(K)\longrightarrow MaxSperH(K)$ est donn\'{e}e par $\alpha
\longmapsto \alpha ^{\max }$ (o\`{u} $\alpha ^{\max }$ est l'unique
sp\'{e}cialisation maximale de $\alpha )$ ;

\bigskip

$res:M(K)\longrightarrow Hom(H(K),\mathbb{R})$ est donn\'{e}e par $\lambda
\mapsto \lambda _{\mid H(K)}$ ;

\bigskip

$j:Hom(H(K),\mathbb{R})\longrightarrow MaxSperH(K)$ est donn\'{e}e par $%
\varphi \mapsto \alpha _{\varphi }$ (o\`{u}, selon la notation adopt\'{e}e
pour le spectre r\'{e}el, $\alpha _{\varphi }=\varphi ^{-1}($ $\mathbb{R}%
^{2}),$ ou $\alpha _{\varphi }=(\ker \varphi ,\overline{\alpha })$ avec $%
\overline{\alpha }=$ $\mathbb{R}^{2}\cap quot(\varphi (H(K)).$

\bigskip

Tous ces espaces sont compacts et les topologies de $M(K)$ et $MaxSperH(K)$
sont les topologies quotients issues de $\Lambda $ et $sp.$

L'espace $\chi (K)$ des ordres d'un corps est donc hom\'{e}omorphe \`{a} $%
MinSperH(K),$ et l'espace $M(K)$ des places r\'{e}elles est lui
hom\'{e}omorphe \`{a} $MaxSperH(K).$

\section{Eventails valu\'{e}s et $\mathbb{R}-$places}

\subsection{{\protect\large Compatibilit\'{e} d'un pr\'{e}ordre et d'une
valuation (par ex.[14])}}

\begin{definition}
Un pr\'{e}ordre $T$ dans un corps $K$ est un sous-ensemble $T\subseteq K$
satisfaisant $\ T+T\subseteq T,$ $T.T\subseteq T$, $0,1\in T$, $-1\notin T$
et $T^{\ast }=T\backslash \{0\}$ est un sous-groupe de $K^{\ast }$ .
\end{definition}

\begin{definition}
On dit que le pr\'{e}ordre $T$ est compatible avec une valuation $v$ si $1+%
\QTR{sl}{m}_{v}\subset T,o\grave{u}$ $m_{v}$ d\'{e}signe l'id\'{e}al maximal
de l'anneau de valuation de $v.$ Alors $T$ induit sur le corps r\'{e}siduel $%
k_{v}$ un pr\'{e}ordre $\overline{T}.$
\end{definition}

\subsection{Les \'{e}ventails valu\'{e}s (Jacob [13])}

\begin{definition}
Un \'{e}ventail $T\subset K$ est un pr\'{e}ordre compatible avec une
valuation $v,$ qui induit sur le corps r\'{e}siduel $k$ un pr\'{e}ordre $%
\overline{T},$ tel que $\overline{T}$ est l'intersection de deux ordres
usuels ou est un ordre usuel ($\overline{T}$ est un \'{e}ventail trivial).
\end{definition}

\begin{definition}
Un\emph{\ }\'{e}ventail valu\'{e} est un pr\'{e}ordre $T$ pour lequel il
existe une valuation r\'{e}elle v, compatible avec le pr\'{e}ordre ($%
1+m_{v}\subset T$), qui induit un ordre archim\'{e}dien sur le corps
r\'{e}siduel $k_{v}$.
\end{definition}

\textbf{Remarque} : on appelle \'{e}ventail trivial l'intersection de deux
ordres usuels ou un ordre usuel.

\textbf{Exemples :}

\QTP{Body Math}
1- Les ordres usuels $P$ sont des \'{e}ventails valu\'{e}s (de niveau $1,$
i.e$.\sum K^{2}\subset P)$.

\QTP{Body Math}
2- Les ordres $P_{n}$ de niveau $n$ sont des \'{e}ventails valu\'{e}s (de
niveau $n)$.

\QTP{Body Math}
3- Soit $\Lambda ^{-1}(\lambda )=\{P_{i}\mid \lambda _{P_{i}}=\lambda \}$
(o\`{u} $\Lambda $ est l'application $:\chi (K)\ \longrightarrow M(K)$
d\'{e}finie par $P\mapsto \lambda _{P}),$ alors $T=\cap P_{i}$est un
\'{e}ventail valu\'{e} de niveau 1 minimal.

\subsection{Lien avec les $\mathbb{R}-$places (Becker-Berr-Gondard [4])}

Comme dans le cas des ordres usuels, on peut associer \`{a} tout
\'{e}ventail valu\'{e} une $\mathbb{R}-$place$.:$

Soit $K$ un corps ordonnable, $T$ un \'{e}ventail valu\'{e} de $K$ $;$ on
sait que 
\begin{equation*}
A(T)=\{a\in K\mid \exists r\in \mathbb{Q},\text{ }r\pm a\in T\}
\end{equation*}
est un anneau de valuation d'id\'{e}al maximal $I(T)=\{a\in K\mid \forall
r\in \mathbb{Q}^{\ast },$ $r\pm a\in T\},$ et que $T$ induit sur le corps r%
\'{e}siduel $k_{v}=A(T)/I(T)$ un ordre archim\'{e}dien $\overline{T}$ . $%
(k_{v,}\overline{T})$ se plonge donc avec unicit\'{e} dans $(\mathbb{R},%
\mathbb{R}^{2})$ $;$ on note $i$ ce plongement et $\pi $ l'application
canonique de $K$ dans $k_{v}\cup \{\infty \}$ (o\`{u} $\pi (a)=\infty )$ si $%
a\notin A(T)).$

\begin{definition}
La $\mathbb{R}-place$ associ\'{e}e \`{a} $T,$ $\lambda _{T}:K\rightarrow 
\mathbb{R}\cup \{\infty \},$ est d\'{e}finie par le diagramme commutatif :
\end{definition}

$\qquad \qquad \qquad \qquad \qquad \qquad \qquad K\qquad \overset{\lambda
_{T}}{\longrightarrow }\qquad \mathbb{R}\cup \{\infty \}$ \ \ \ 

\bigskip

$\qquad \qquad \qquad \qquad \qquad \qquad \qquad \pi \searrow \qquad \qquad
\nearrow i$

\qquad \qquad \qquad \qquad \qquad \qquad \qquad \qquad $k_{v}\cup \{\infty
\}$

\subsection{Cl\^{o}tures pour une $\mathbb{R-}$place}

Les \emph{corps henseliens r\'{e}siduellement r\'{e}els-clos} ont
\'{e}t\'{e} pr\'{e}sent\'{e}s par Becker-Berr-Gondard dans [4] comme des
cl\^{o}tures par extension alg\'{e}brique d'un corps muni d'un \'{e}ventail
valu\'{e} ; on peut aussi utiliser cette approche pour obtenir des
cl\^{o}tures par extensions alg\'{e}briques pour une $\mathbb{R-}$place
donn\'{e}e $\lambda $ en lui associant l'\'{e}ventail valu\'{e} minimal de
niveau 1 qui lui est associ\'{e} : $T=\cap P_{i},$ o\`{u} les $P_{i}$ sont
dans la fibre $\Lambda ^{-1}(\lambda )=\{P_{i}\in \chi (K)\mid \lambda
_{P_{i}}=\lambda \}$.

\bigskip

Tous les corps henseliens r\'{e}siduellement r\'{e}els-clos sont clos pour
leur unique $\mathbb{R-}place.$

C'est aussi vrai pour les cas particuliers que sont les \emph{corps
r\'{e}els-clos g\'{e}n\'{e}ralis\'{e}s} de [5] et les \emph{corps de Rolle}
de [9].

\bigskip

C'est seulement dans le cas o\`{u} $\Lambda ^{-1}(\lambda )=\{P\in \chi
(K)\mid \lambda _{P}=\lambda \}$ ne contient qu'un seul ordre que l'on
obtient une cl\^{o}ture pour la $\mathbb{R-}$place unique \`{a} isomorphisme
pr\`{e}s (cas des corps r\'{e}els-clos). Dans tous les autres cas, les
classes d'isomorphisme des cl\^{o}tures correspondent au choix d'une
cha\^{i}ne infinie d'\'{e}ventails valu\'{e}s associ\'{e}e \`{a} la $\mathbb{%
R-}$place (voir [4])$.$

\subsection{Les signatures g\'{e}n\'{e}ralis\'{e}es}

La notion d'\'{e}ventail valu\'{e} permet de d\'{e}finir celle de signature
g\'{e}n\'{e}ralis\'{e}e.

\begin{definition}
(Schwartz [20]) : Une signature g\'{e}n\'{e}ralis\'{e}e est un morphisme de
groupes ab\'{e}liens $\sigma :K^{\ast }\longrightarrow G$ dont le noyau est
tel que $T=\ker \sigma \cup \left\{ 0\right\} $ est un \'{e}ventail
valu\'{e}.
\end{definition}

\textbf{Exemples} :

1- Si $\sigma ${\Large \ }est un morphisme de groupes, $\sigma $ : $K^{\ast
}\,\,\longrightarrow \left\{ \pm 1\right\} $ de noyau additivement
ferm\'{e}, alors $\sigma $ est une signature et $P=ker$ $\sigma \cup \left\{
0\right\} $ est un ordre.

2- Si $\sigma :K^{\ast }\rightarrow \mu _{2n}$ morphisme de groupes
ab\'{e}liens de noyau additivement ferm\'{e}, alors $P=\ker \sigma \cup
\{0\} $ est un ordre dont le niveau divise $n$.

\section{Espace d'ordres abstrait et places abstraites}

L'espace des ordres d'un corps - \'{e}tudi\'{e} en relation avec les formes
quadratiques et les valuations r\'{e}elles - a \'{e}t\'{e} \`{a} l'origine
de la th\'{e}orie des espaces d'ordres abstraits (M. Marshall 1979-80).

\subsection{\protect\large Espaces d'ordres abstraits }

\begin{definition}
(Marshall [15]) : Un espace d'ordres abstrait est un couple $(X,G)$ o\`{u} G
est un groupe d'exposant $2$ (donc ab\'{e}lien), $-1$ un \'{e}l\'{e}ment
distingu\'{e} de $G,$ et $X$ \ un sous-ensemble de $Hom(G,\{1,-1\})$
satisfaisant les quatre axiomes suivants :

(1) $X$ est un sous-ensemble ferm\'{e} de $Hom(G,\{1,-1\})$

(2) $\forall \sigma \in X$ \ $\sigma (-1)=-1$

(3) $\underset{\sigma \in X}{\cap }\ker \sigma =\left\{ 1\right\} $ \
(o\`{u} $\ker \sigma =\{a\in G\mid \sigma (a)=1\}$)

(4) Pour des formes quadratiques $f,g$ sur $G,$%
\begin{equation*}
D_{X}(f\oplus g)=\cup \{D_{X}\left\langle x,y\right\rangle \mid x\in
D_{X}(f),y\in D_{X}(g)\}
\end{equation*}
o\`{u}\textit{\ }$D_{X}(f)=\{a\in G$\textit{\ }$\mid a$ \textit{%
repr\'{e}sent\'{e} par }$f\},$\textit{\ i.e. il existe }$g$\textit{\ telle
que }$f\equiv _{X}\left\langle a\right\rangle \oplus g$\textit{\ (o\`{u} }$%
f\equiv _{X}h$\textit{\ ssi }$f$\textit{\ et }$h$\textit{\ ont m\^{e}me
dimension, et }$\forall \sigma \in X$ \textit{m\^{e}me signature )}
\end{definition}

\QTP{Body Math}
Si on consid\`{e}re les signatures, un \'{e}ventail de niveau 1 \`{a} quatre
\'{e}lements est caract\'{e}ris\'{e} par $\sigma _{0}\sigma _{1}\sigma
_{2}\sigma _{3}=1$

\QTP{Body Math}
Dans le cadre abstrait\textbf{\ }un \'{e}ventail est un espace d'ordre
abstrait $(X,G)$ tel que $X=\left\{ \sigma \in \widehat{G}\mid \sigma
(-1)=-1\right\} .$

\QTP{Body Math}
Il est caract\'{e}ris\'{e} par : $\forall \sigma _{0},\sigma _{1},\sigma
_{2}\in X$ on a $\sigma _{0}\sigma _{1}\sigma _{2}\in X)$

\subsection{{\protect\large Espaces de signatures abstraits (cas }$2^{n})$}

\begin{definition}
\textbf{\ }Un espace de signatures abstrait est un couple $(X,G)$, $G$
groupe ab\'{e}lien d'exposant $2^{n}$, $X$ $\subset Hom(G,\mu _{2^{n}})=\chi
(G)$ tel que :

(0) $\forall \sigma \in X,$ $\forall k\in \mathbb{N}$ avec $k$ impair, $%
\sigma ^{k}\in X$ ;

(1) $X$ est un sous-ensemble ferm\'{e} de $\chi (G)$ ;

(2) $\forall \sigma \in X$ \ $\sigma (-1)=-1$ ($-1$ \'{e}l\'{e}ment
distingu\'{e}$)$ ;

(3) $\underset{\sigma \in X}{\cap }\ker \sigma =\left\{ 1\right\} $\ (o\`{u} 
$\ker \sigma =\{a\in G\mid \sigma (a)=1\}$)

(4) des formes $f,g$ sur $G,$%
\begin{equation*}
D_{X}(f\oplus g)=\cup \{D_{X}\left\langle x,y\right\rangle \mid x\in
D_{X}(f),y\in D_{X}(g)\}
\end{equation*}
\end{definition}

\subsection{\protect\large P-Structures}

\begin{definition}
(Marshall [15]): Une\emph{\ }P-structure est une relation d'\'{e}quivalence
sur un espace de signatures $(X,G)$ telle que l'application canonique $%
\Lambda :X\rightarrow M$ (o\`{u} $M$ est l'ensemble des classes
d'\'{e}quivalence) v\'{e}rifie

(1) chaque fibre est un \'{e}ventail ;

(2) si $\sigma _{0}\sigma _{1}\sigma _{2}\sigma _{3}=1$ alors $\left\{
\sigma _{0},\sigma _{1},\sigma _{2},\sigma _{3}\right\} $ a une intersection
non vide avec au plus deux fibres.
\end{definition}

Tout espace d'ordres abstrait poss\`{e}de une $P$-structure $M$ (voir [15])$%
. $ Mais $M$ muni de la topologie quotient n'est pas toujours un espace
s\'{e}par\'{e}.

\subsection{Places r\'{e}elles abstraites}

Nous avons montr\'{e} que dans certaines conditions on peut associer \`{a}
l'espace d'ordres abstrait une ``P-structure`` correspondant \`{a} l'espace
des $\mathbb{R}-$places dans le cas des corps.

\begin{theorem}
(Gondard-Marshall [13]) : Soit $(X,G)$ sous espace d'un espace de signatures 
$(X^{\prime },G^{\prime })$ d'exposant $2^{n}$, $avec$ $n\geq 1$.

Pour $\sigma ,\tau \in X$ on dit que \ $\sigma \sim \tau $ si $\sigma \tau
\in X^{\prime 2}.$ Alors sont \'{e}quivalents :

(1) $\sim $ d\'{e}finit une $P$-structure sur $X.$

(2) Si $\sigma _{0}\sigma _{1}\sigma _{2}\sigma _{3}=1$, alors soit $\sigma
_{0}$ est li\'{e} par $\sim $ \`{a} exactement un des $\sigma _{1},\sigma
_{2},\sigma _{3},$ soit $\sigma _{0}$ est li\'{e} par $\sim $ \`{a} tous les 
$\sigma _{1},\sigma _{2},\sigma _{3}.$

De plus dans ce cas la $P$-structure d\'{e}finie sur $X$ par $\sim $ a une
topologie s\'{e}par\'{e}e.
\end{theorem}

\textbf{Origine du r\'{e}sultat}

Dans le cas des corps, l'espace des $\mathbb{R}-$places est connu d\`{e}s
qu'on connait les ordres usuels et les ordres de niveau 2.

Ou, de mani\`{e}re \'{e}quivalente, d\`{e}s qu'on connait les \'{e}ventails
valu\'{e}s de niveau 1 minimaux $T=\cap \{P\in \chi (K)\mid \lambda
_{P_{{}}}=\lambda \}$

Par exemple $\mathbb{Q(}\sqrt{2})$ et $\mathbb{R((}X\mathbb{))}$ ont des
espaces d'ordres isomorphes, mais le premier a deux $\mathbb{R}-$places et
pas d'ordres de niveau $2,$ alors que le second a une seule $\mathbb{R}-$%
place et une cha\^{i}ne d'ordres de niveaux puissances de $2.$

Les cha\^{i}nes d'ordres de niveau puissance de 2 commencent par une paire
d'ordres usuels $P_{0},P_{1}$ et l'ordre $P_{2}$ de niveau $2$\ qui leur
correspond satisfait $P_{2}\cup -P_{2}=(P_{0}\cap P_{1})\cup (-(P_{0}\cap
P_{1})),$ d'o\`{u} on peut d\'{e}duire, puisque $a^{2}\in
P_{2}\Longleftrightarrow a\in \pm P_{2\text{ }}$, que $\sigma
_{2}^{2}(a)=\sigma _{2}(a^{2})=1\Longleftrightarrow \sigma _{0}(a)\sigma
_{1}(a)=1$ (en notant $P_{i}=\ker \sigma _{i}\cup \left\{ 0\right\} ).$

\section{Quelques questions ouvertes}

\subsection{Espaces des $\mathbb{R-}$places}

On sait peu de choses sur l'espace des $\mathbb{R-}$places \`{a} part le
fait qu'il est compact s\'{e}par\'{e}. Quels espaces topologiques compacts
s\'{e}par\'{e}s peuvent \^{e}tre des espaces de $\mathbb{R}-$places\ ?

D'autre part si on connait l'espace des $\mathbb{R-}$places d$^{\prime }$un
corps $K,$ que peut on dire de l'espace des $\mathbb{R-}$places des corps
extensions de $K$ ?

D'apr\`{e}s les travaux de Schulting [17] on sait que pour $%
L=K(X_{1},...,X_{n})$ et pour $L=K((X)),$ la connexit\'{e} \'{e}ventuelle de 
$M(K)$ passe \`{a} $M(L),$ et r\'{e}ciproquement ; on sait aussi que si tous
les ordres de $K$ s'\'{e}tendent \`{a} $L$ et si $M(L)$ est connexe, alors $%
M(K)$ est connexe.

On ne sait rien sur le cas des extensions alg\'{e}briques en g\'{e}n\'{e}ral.

\subsection{$\mathbb{R}-places${\protect\large \ g\'{e}om\'{e}triques}}

Dans le cadre g\'{e}om\'{e}trique sur $\mathbb{R},$ $K=\mathbb{R}(Y)$ , $Y$
vari\'{e}t\'{e} de dimension $d$ , on appellera \textit{\'{e}ventail
valu\'{e} g\'{e}om\'{e}trique minimal de niveau 1, } un \'{e}ventail
valu\'{e} $T=\cap P_{i}$ o\`{u} il y a $2^{d}$ ordres $P_{i}$ de $\mathbb{R}%
(Y)$ distincts, et tel que pour $v$ la plus fine des valuations compatibles
avec $T$ (la valuation de Becker dont l$^{\prime }$anneau de valuation est $%
A(T)$)$,$ $T$ induit sur le corps r\'{e}siduel un ordre archim\'{e}dien $%
\overline{T}$.

Avec F. Acquistapace et F. Broglia, nous avons consid\'{e}r\'{e} les $%
\mathbb{R}-$places g\'{e}om\'{e}triques correspondant \`{a} ces
\'{e}ventails valu\'{e}s de niveau 1 minimaux particuliers (qui sont aussi
certains des \'{e}ventails g\'{e}om\'{e}triques d\'{e}finis par Andradas et
Ruiz [1]).

Dans le cas des courbes projectives lisses il est connu que $V(\mathbb{R)}$,
l'ensemble des points r\'{e}els, est hom\'{e}omorphe \`{a} $M(K)$ (dans ce
cas toutes les places r\'{e}elles sont g\'{e}om\'{e}triques).

Dans le cas g\'{e}n\'{e}ral, la fibre des $\mathbb{R-}$places
g\'{e}om\'{e}triques centr\'{e}es sur un point r\'{e}gulier de la
vari\'{e}t\'{e} doit pouvoir \^{e}tre interpr\'{e}t\'{e}e
g\'{e}om\'{e}triquement et d\'{e}crire le voisinage du point dans la
vari\'{e}t\'{e}.

\subsection{Espaces des \'{e}ventails valu\'{e}s}

Il existe une bijection entre l'espace des $\mathbb{R-}$places et les
\'{e}ventails valu\'{e}s minimaux (minimaux comme pr\'{e}ordres) de niveau 1.

A tout \'{e}ventail valu\'{e} on associe une et une seule $\mathbb{R-}$place$%
.$

\bigskip

R\'{e}ciproquement, on a l'application $\Lambda :\chi (K)\ \longrightarrow
M(K)$ ,\ donn\'{e}e par $P\mapsto \lambda _{P}$ $;$ si on consid\`{e}re
l'image r\'{e}ciproque $\Lambda ^{-1}(\lambda )=\{P_{i}\in \chi (K)\mid
\lambda _{P_{i}}=\lambda \}$ alors $T=\cap P_{i}$ est un \'{e}ventail
valu\'{e} de niveau 1 minimal.

Les ordres peuvent aussi \^{e}tre consid\'{e}r\'{e}s comme des \'{e}ventails
valu\'{e}s de niveau 1 maximaux (comme pr\'{e}ordres).

\bigskip

D'apr\`{e}s le th\'{e}or\`{e}me 18 on peut penser \`{a} associer \`{a}
l'espace des \'{e}ventails valu\'{e}s de niveau 1 le spectre r\'{e}el de
l'anneau d'holomorphie r\'{e}el $Sper(H(K))$.

Ces consid\'{e}rations pourraient sans doute permettre de d\'{e}finir une
notion d'espace d'\'{e}ventails valu\'{e}s abstrait.

\subsection{Probl\`{e}me de Marshall}

Rappelons l'\'{e}nonc\'{e} du probl\`{e}me de Marshall :

`` Tout espace d'ordre abstrait est\/-il l'espace des ordres d'un corps?''.

Pouvoir obtenir une th\'{e}orie abstraite satisfaisante des espaces de
places r\'{e}elles peut \^{e}tre un premier pas vers sa r\'{e}solution. En
effet si un espace d'ordres abstrait est r\'{e}alis\'{e} comme espace des
ordres d'un corps, on devra pouvoir trouver un espace abstrait de places
r\'{e}elles correspondant aux $\mathbb{R-}$places du corps. Une variante
\`{a} cette id\'{e}e serait de chercher une notion d'espace d'\'{e}ventails
valu\'{e}s abstrait comme sous famille de l'espace des \'{e}ventails
abstraits d'un espace d'ordres abstrait.

\bigskip

\textbf{R\'{e}f\'{e}rences}

\bigskip

[1] C. Andradas et J. Ruiz : \textit{Algebraic and analytic geometry of fans,%
} Memoirs of AMS, vol.115 \#553, (1995).

[2] E. Becker : \textit{Valuations and real places in the theory of formally
real fields }\textbf{et }\textit{The real holomorphy ring and sums of }$%
2^{n} $\textit{-th powers}, in LNM\ 959, G\'{e}om\'{e}trie R\'{e}elle et
Formes Quadratiques, 1-40 \textbf{et} 139-181, (1982).

[3] E. Becker : \textit{Extended Artin-Schreier theory of fields, }Rocky
Mountain J. of Math., vol 14\#4, 881-897 (1984).

[4] E. Becker, R. Berr, et D. Gondard : \textit{Valuation fans and
residually real-closed henselian fields, }J. Algebra,215, 574-602 (1999).

[5] E. Becker, R. Berr, F. Delon, et D. Gondard : \textit{Hilbert's 17th
problem for sums of 2n-th powers, }J. Reine Angew. Math. 450, 139-157 (1994).

[6] E. Becker et D. Gondard : \textit{Notes on the space of real places of a
formally real field, }in RAAG (Trento), W. de Gruyter\textit{,} 21-46,
(1995).

[7] J.--L. Colliot Th\'{e}l\`{e}ne : \textit{Eine Bemerkung zu einem Satz
von E. Becker und D. Gondard}, http://www.math.psud.fr/\symbol{126}colliot/
; pr\'{e}publication \#5, \`{a} para\^{i}tre dans Math. Zeitschrift.

[8] T. Craven :\textit{\ The Boolean space of orderings of a field, }Trans.
Amer. Math.Soc. 209 , 225-235 (1975).

[9] D. Gondard-Cozette : \textit{Axiomatisations simples des th\'{e}ories
des corps de Rolle, }Manuscripta Mathematica 69, 267-274 (1990).

[10] D. Gondard et M. Marshall : \textit{Towards an abstract description of
the space of real places, }Contemporary Mathematics, vol 253, AMS, 79-113,
(2000).

[11] D. Gondard et M. Marshall : \textit{Real holomorphy rings and the
complete real spectrum, }en pr\'{e}paration.

[12] J. Harman : \textit{Chains of higher level orderings, }Contemporary
Mathematics 8, 141-174, (1982).

[13] B. Jacob : \textit{Fans, real valuations, and hereditarily-pythagorean
fields, }Pacific J. Math. 93, 95-105 (1981).

[14] T. Y. Lam : \textit{Orderings,} \textit{Valuations and Quadratic Forms,
AMS, Regional Conference Series in Mathematics \#52, (1983).}

[15] M. Marshall : \textit{Spaces of Orderings and Abstract Real Spectra, }%
LNM\textit{\ \ 1636, }Springer\textit{-}Verlag\textit{, (1996).}

[16] C. Scheiderer : \textit{A short remark on a theorem by Becker and
Gondard,} http://www.uni-duisburg.de/FB11/FGS/F1/claus.html\#notes

[17] H.-W. Sch\"{u}lting :\textit{\ On real places of a field and their
holomorphy ring , }Communications in Algebra, 10, 1239-1284 (1982).

[18] H.-W. Sch\"{u}lting : \textit{The strong topology on real algebraic
varieties, }Contemporary Mathematics 8, 141-174, (1982).

[19] H.-W. Sch\"{u}lting :\textit{\ Real holomorphy rings in real algebraic
geometry, }LNM\ 959, G\'{e}om\'{e}trie R\'{e}elle et Formes Quadratiques,
1-40 \textbf{et} 139-181, (1982).

[20] N. Schwartz : \textit{Signatures and real closures of fields, }in
S\'{e}minaire Structures Alg\'{e}briques Ordonn\'{e}es, Publications de
l'Universit\'{e} Paris 7, vol 33, 65-78 (1990).

\bigskip

D. Gondard-Cozette

Institut de Math\'{e}matiques de Jussieu (UMR 7586 du CNRS)

Sorbonne Universit\'{e}

4, Place Jussieu

75252 Paris cedex 05 (France)

\end{document}